\documentclass[12pt,a4paper]{amsart}
\usepackage{amssymb,amsmath,amsthm}

\usepackage{graphicx}

\usepackage{epsfig}

\long\def\onefigure#1#2{
\begin{figure*}[tbp]
\begin{center}
#1
\end{center}
\caption{#2}
\end{figure*}
} 

\newcommand{\lipefig}[2]  
{\onefigure{\mbox{\psfig{file=#1.eps}}}{\label{f:#1} #2} }

\begin{document}

\theoremstyle{plain}
\newtheorem{theorem}{Theorem}[section]
\newtheorem{lemma}[theorem]{Lemma}
\newtheorem{prop}[theorem]{Proposition}
\newtheorem{corollary}[theorem]{Corollary}
\newtheorem{claim}[theorem]{Claim}

\newcommand{\sgn}{\textrm{sign}}
\newcommand{\tri}{\triangle}
\newcommand{\al}{\alpha}
\newcommand{\be}{\beta}
\newcommand{\ga}{\gamma}
\newcommand{\la}{\lambda}
\newcommand{\eps}{\varepsilon}
\newcommand{\N}{\mathbb{N}}
\newcommand{\z}{{\mathbb{Z}^2}}
\newcommand{\R}{\mathbb{R}}
\newcommand{\A}{\mathrm{A}}
\newcommand{\cc}{\mathcal{C}}
\newcommand{\xx}{\mathcal{A}}
\newcommand{\pp}{\mathcal{G}}
\newcommand{\dd}{\mathcal{D}}
\newcommand{\hh}{\mathcal{H}}
\newcommand{\dist}{\textrm{dist}}
\newcommand{\conv}{\textrm{conv}}
\newcommand{\intt}{\textrm{int}}
\newcommand{\sign}{\textrm{sign}}

\numberwithin{equation}{section}

\title{Block partitions of sequences}

\author{Imre B\'ar\'any, Victor S. Grinberg}



\begin{abstract} Given a sequence $A=(a_1,\dots,a_n)$ of real numbers, a block $B$ of the $A$ is either a set $B=\{a_i,a_{i+1},\dots,a_j\}$  where $i\le j$ or the empty set. The size $b$ of a block $B$ is the sum of its elements. We show that when each $a_i \in [0,1]$ and $k$ is a positive integer, there is a partition of $A$ into $k$ blocks $B_1,\dots,B_k$ with $|b_i-b_j|\le 1$ for every $i,j$. We extend this result in several directions.
\end{abstract}

\maketitle
\bigskip

\section{Introduction}
Assume $A=(a_1,\dots,a_n)$ is a sequence of real numbers $a_i \in [0,1]$. A block $B$ of the sequence is either the empty set or it is  $\{a_i,a_{i+1},\dots,a_j\}$  with $i\le j$. The size of the block $B$, to be denoted by $b$, is just the sum of the elements in $B$. Blocks $B_1,\dots,B_k$ form a partition of $A$ if every element of $A$ belongs to exactly one block. We always assume that if $a_h$ is the last element of a non-empty block, then $a_{h+1}$ is the first element of the next non-empty block.

It is easy to see that, for a given $k\in \N$, there is a $k$-partition of $A$ into blocks $B_1,\dots,B_k$, of sizes $b_1,\dots,b_k$, such that

\begin{equation}\label{eq:Mm}
|b_i-b_j|\le 2 \mbox{ for all }i,j \in [k].
\end{equation}

Here and later, $[k]$ stands for the set $\{1,2,\dots,k\}$. To see this define $S_j=\sum_1^ja_i$ and set $S=S_n$. The condition $a_i \in [0,1]$ implies that for every $h \in [k-1]$ there is a subscript $m(h)$ such that $hS/k-1/2\le S_{m(h)} \le hS/k+1/2$ and $m(h)$ is a non-decreasing function of $h$. The partial sums $S_0,S_{m(1)},\dots,S_{m(k-1)},S_n$ split $A$ into $k$ blocks $B_1,\dots,B_k$ that satisfy $S/k-1\le b_i\le S/k+1$ for all $i$ and consequently $|b_i-b_j|\le 2$ for all $i,j$.

In the first part of this paper we show the existence of a $k$-partition with $|b_i-b_j|\le 1$. Then we extend this result to infinite sequences. Finally we show that the bound in (\ref{eq:Mm}) holds under much weaker conditions. Related problems are treated in \cite{BD} and \cite{BG}.

\section{Finite sequences}

Our starting result is the following theorem.

\begin{theorem}\label{th:main} Given a sequence  $A=(a_1,\dots,a_n)$ of real numbers $a_i \in [0,1]$ and a positive integer $k$, there is a partition of $A$ into $k$ blocks $B_1,\dots,B_k$ with $\max b_i \le \min b_i +1$.
\end{theorem}

{\bf Remarks.} This result is best possible in the sense that, in general, $\max b_i - \min b_i$ cannot be made smaller than 1. There are many examples showing this, for instance $A=(1/2,1,\dots,1,1/2)$ with $k-1$ ones in the middle, or when every $a_i=1$ and $k$ does not divide $n$. When $k>n$, the last example shows also that empty blocks have to be allowed. But no empty block can be present when $\sum_1^na_i>k$, and actually even when $\sum_1^na_i>k-1$.

\bigskip
{\bf Proof.} Given a $k$-partition $P$ of $A$ with blocks $(B_1,\dots,B_k)$ let $M(P)=\max_{i \in [k]} b_i$ and $m(P)=\min_{i \in [k]} b_i$. We give an algorithm that finds the required partition. It starts with an arbitrary $k$-partition $P$. On each iteration, the current partition $P$ is changed to another one, $P^*$, and the only difference is that either the last element of $B_h$ is moved to $B_{h+1}$ or the first element of $B_h$ is moved to $B_{h-1}$ for a unique $h \in [k]$.

Here comes the algorithm plus some comments.

\medskip
\begin{enumerate}
\item[(1)] Fix $p \in [k]$ with $M(P)=b_p$. So $B_p$ is a maximal block of $P$.

\item[(2)] If $M(P)\le m(P)+1$, then stop.

\item[(3)] If $M(P)> m(P)+1$, then let $m(P)=b_q$ for the $q \in [k]$ which is closest to $p$ (ties broken arbitrarily). Thus $B_q$ is a minimal block of $P$. Let $B_h$ be the block next to $B_q$ between $B_p$ and $B_q$. (Note that $B_h$ is a non-empty block: if it were, then $m(P)=0$ and we should have chosen $B_h$ instead of $B_q$.) So either $p< q$ and then $h=q-1$ and we define $P^*$ by moving the last element from $B_h=B_{q-1}$ to $B_q$, or $q< p$, and then $h=q+1$ and $P^*$ is obtained by moving the first element of $B_h=B_{q+1}$ to $B_q$. Set $P=P^*$. If $p=h$, then go to (1), else go to (2).
\end{enumerate}

\medskip
We prove next that this algorithm terminates with the required partition. Note first that in step (3) the size of every block in $P^*$ is at most $M(P)$. Indeed, the size of $B_h$ does not increase, and the size of $B_q$ increases by some $a_i\le 1$ and since we have $m(P)+1<M(P)$, $m(P)+a_i<M(P)$ follows. This shows that $M(P^*)\le M(P)$, that is, the size of maximal block does not increase during the algorithm. Note also that no new block of size $M(P)$ is created in step (3).

\begin{claim} Step (3) is repeated at most $kn$ times with $B_p$ being the same block in $P^*$ and in $P$.
\end{claim}

For the {\bf proof}, let us define $f(P,p)=\sum_{i=1}^k |i-p||B_i|$ where, as usual, $|B_i|$ denotes the number of elements in $B_i$. It is clear that $f(P,p)$ takes positive integral values and is always less than $kn$. It is also evident that $f(P,p)<f(P^*,p)$, which proves the claim.\hfill$\Box$

\medskip
Thus after at most $kn$ iteration of (3), the algorithm decreases the size of $B_p$ and so goes to (1). Consequently it decreases either the number of maximal blocks or $M(P)$. As there are only finitely many block partitions, the algorithm eventually terminates with (2).\hfill$\Box$

\medskip
\begin{prop} The above algorithm takes at most $O(kn^3)$ steps.
\end{prop}

The {\bf proof} follows from three simple facts. Note that a block is just a set of consecutive elements of the sequence.
\begin{enumerate}
\item[(a)] No block can be fixed at step (1) more than once.
\item[(b)] No block can serve as maximal block for more than $kn$ iterations of the loop ``go to (2)'' (and due to (a), in total, as well).
\item[(c)] There are no more than $O(n^2)$ blocks.
\end{enumerate} \hfill$\Box$

\bigskip
Theorem \ref{th:main} can be strengthened by removing the condition $a_i\ge 0$:

\begin{theorem}\label{th:negative} Given a sequence  $A=(a_1,\dots,a_n)$ of real numbers $a_i \le1$ for every $i \in [n]$ with $S=\sum_1^na_i\ge 0$ and a positive integer $k$, there is a partition of $A$ into $k$ blocks $B_1,\dots,B_k$ with $\max b_i \le \min b_i +1$.
\end{theorem}

The proof is based on a lemma that can suitably preprocess the sequence $A$.

\begin{lemma}\label{l:preproc} Given a sequence  $A=(a_1,\dots,a_n)$ of real numbers $a_i \le1$ for every $i \in [n]$ with $\sum_1^na_i\ge 0$, there is a partition of $A$ into blocks $(C_1,\dots,C_m)$ such that $c_i=\sum_{a_j\in C_i}a_j\in [0,1]$ for every $i \in [m]$.
\end{lemma}

The {\bf proof} is by induction on $n$ and the case $n=1$ is trivial. Assume the statement holds for sequences with fewer than $n$ entries ($n\ge 2$). We show that the statement holds for  $A=(a_1,\dots,a_n)$. If $S\le 1$, then we can choose a single block $C_1=A$. Otherwise $S>1$ and we choose the smallest subscript $h\in [n]$ such that the size, $c_1$, of the block $C_1=(a_1,\dots,a_h)$ is positive. It is clear that $c_1 \in (0,1]$. The sequence $A^*=(a_{h+1},\dots,a_n)$ has fewer than $n$ elements, every $a_i\le 1$, and the sum of its elements is $S-c_1>0$. So induction gives a partition of $A^*$ into blocks $(C_2,\dots,C_m)$ with all $c_i \in [0,1]$. They, together with $C_1$ form the required partition of $A$.\hfill$\Box$

\bigskip
{\bf Remark.} There is a simple algorithm that produces the partition $(C_1,\dots,C_m)$. Namely, starting with $A=(a_1,\dots,a_n)$, check if there is an $i \in [n-1]$ with $a_ia_{i+1}\le 0$. If there is no such $i$, then the partition with blocks $C_i=(a_i)$ satisfies the requirements. If there is such an $i$ replace $A$ by the sequence $(a_1,\dots,a_{i-1},a_i+a_{i+1},a_{i+2},\dots,a_n)$ of length $n-1$ and continue. The algorithm terminates either with the sequence $(0)$ consisting a single zero, or with a sequence $(c_1,\dots,c_m)$ where each $c_i \in (0,1]$. Note that preprocessing takes $O(n)$ iterations with this algorithm.

\medskip
The {\bf proof} of Theorem~\ref{th:negative} is quite easy now. Just apply the preprocessing lemma to $A$ to obtain the partition into blocks $(C_1,\dots,C_m)$. The sequence $C=(c_1,\dots,c_m)$ satisfies the conditions $c_i \in [0,1]$ so Theorem \ref{th:main} applies and gives the suitable partition of $C$ which is, in fact, a suitable partition of $A$ as well.\hfill$\Box$

\begin{corollary} Given a sequence $A = (a_1,...,a_n)$ of real numbers with
$a_i\in[-1, 1]$ for all $i$ and a positive integer $k$, there is a partition of $A$ into $k$
blocks $B_1,...,B_k$ such that $\max b_i - \min b_i \le 1$.
\end{corollary}

The {\bf proof} follows immediately from Theorem~\ref{th:negative} if $a_1+...+a_n\ge 0$. When $a_1+...+a_n < 0$, replace each $a_i$ by $-a_i$, and apply the same theorem. The resulting block partition is a block partition of the original sequence which satisfies $\max b_i - \min b_i \le 1$. \hfill$\Box$

\section{Infinite sequences}

Assume now that $A=(a_1,a_2,\dots)$ is an infinite sequence of real numbers $a_n \in [0,1]$ and $a$ is a non-negative real. To extend our main theorem, we wish to find a partition of $A$ into blocks $(B_1,B_2,\dots)$ such that $\inf b_n \le a \le \sup b_n \le \inf b_n+1$. This may not be possible if $\sum a_n$ is finite: for instance with $\sum a_n=1000$ and $a=400$ the size of the blocks must lie in [399,401]. No set of blocks of this type can partition $A$, clearly. The case is different when $\sum a_n=\infty$.

\begin{theorem}\label{th:infty} Given a sequence  $A=(a_1,a_2,\dots)$ of real numbers $a_i \in [0,1]$ with $\sum a_n=\infty$ and a real number $a\ge 0$, there is a partition of $A$ into blocks $B_1,B_2,\dots$ with $\inf b_i \le a \le \sup b_i \le \inf b_i+1$.
\end{theorem}

{\bf Proof.} The case $a=0$ is easy: just choose $B_1$ to be the empty block and $B_i=(a_{i-1})$, a singleton, for $i=2,3,\dots$. So assume $a>0$. Recall that $S_n=\sum_1^n a_j$.

For every $k\in \N$ let $n(k)$ be the smallest subscript with
\begin{equation}\label{choice}
ka \le S_{n(k)} <ka+1, \mbox{ so } S_{n(k)}=ka+\eps(k),
\end{equation}
where $\eps(k)\in [0,1)$. We can apply Theorem \ref{th:main} to the finite sequence $A_k=(a_1,\dots,a_{n(k)})$ giving a $k$-partition of $A_k$ into blocks $(B_1^k,\dots,B_k^k)$ satisfying
\[
\min_{i \in [k]} b_i^k \le a+\frac {\eps(k)}k \le \max_{i\in [k]} b_i^k \le \min_{i \in [k]} b_i^k+1,
\]
where the middle inequality expresses the fact that the average is between the maximum and the minimum.

Assume now that, for some $k \in \N$, $\min_{i \in [k]} b_i^k \le a$. Then we can produce the required partition as $(B_1^k,\dots,B_k^k,B_{k+1}\dots)$ by defining $B_n$ for $n>k$ recursively as follows. If $B_n$ has been constructed, then let $B_{n+1}$ be the next block with $b_{n+1}\le \min_{i \in [k]} b_i^k+1$.

Assume now that $\min_{i \in [k]} b_i^k>a$ for every $k\in \N$. Then $b_i^k=a+\eps_i^k$ for all $i \in [k]$ and $\sum_1^k \eps_i^k =\eps(k)<1$. This implies that $\max_{i \in [k]} b_i^k =a+ \eps_i^k<a+1$ (for some suitable $i\in [k]$).

It follows that there is an infinite subset $I_1$ of $\N$ such that $B_1^k$ is the same block for all $k \in I_1$. Call this block $B_1$. Further, there is an infinite $I_2\subset I_1$ such that $B_2^k$ is the same block for all $k \in I_2$, call this block $B_2$, etc. We get a partition $(B_1,B_2,\dots)$. Here $\sup b_i\le a+1$ follows from the inequality at the end of the last paragraph.

We show finally that $\inf b_j=a$. Observe first that $b_j=a+\eps_j^k$ for all $k\in I_j$ so we may write $b_j=a+\eps_j$. Then, for $k\in I_j$,
\[
\eps_1+\dots+\eps_j=\eps_1^k+\dots+\eps_j^k \le \eps(k)<1,
\]
showing that $\lim \eps_j=0$. This proves that, indeed, $\inf b_j=a$.\hfill$\Box$

\bigskip
{\bf Remark.} We could have chosen, instead of (\ref{choice}), $n(k)$ as the minimal subscript with
\[
ka-1 < S_{n(k)} \le ka, \mbox{ so } S_{n(k)}=ka-\eps(k),
\]
starting only for $k>1/a$, say. Essentially the same proof works with this choice.

\medskip
Again one can get rid of the condition $a_n\ge 0$.

\begin{theorem}\label{th:inftyneg} Let $A=(a_1,a_2,\dots)$ be a sequence of real numbers $a_i \le 1$ satisfying the condition that for every $k \in \N$ there is $S_n>k$. Then for every real number $a\ge 0$, there is a partition of $A$ into blocks $B_1,B_2,\dots$ with $\inf b_i \le a \le \sup b_i \le \inf b_i+1$.
\end{theorem}

{\bf Proof.} We prove the theorem by reducing it to Theorem~\ref{th:infty}. We construct
a block partition of $A$ so that the size of every block lies in $(0,1]$.
The construction is straightforward.  Find the smallest $i$ such that
$S_i>0$.  Such an $i$ exists because the sequence $S_n$ is not bounded from above.
Clearly, $S_i\le 1$; let $(a_1,\dots,a_i)$ be the first block.  The sequence
$S_n-S_i$ is also unbounded from above, and we
continue this process. Apply Theorem~\ref{th:infty} to the constructed sequence. It is clear that its block partition is in fact a block partition of the sequence $A$ satisfying the requirements. \hfill$\Box$

\section{More general settings}

Assume now that our sequence is $A=(a_1,\dots,a_n)$. Let $s$ be a function defined on the blocks that satisfies the conditions
\begin{enumerate}
\item[(i)] $s(\emptyset)=0$ and $s(B)\ge 0$ for every block $B$,
\item[(ii)] $s(B_1)\le s(B_2)\le s(B_1)+1$ if $B_1\subset B_2$ and $B_2\setminus B_1$ is a singleton.
\end{enumerate}

\begin{theorem}\label{th:gen} Assume $A=(a_1,\dots,a_n)$ and $k$ is positive integer. Under the above conditions on $s$ there is a partition of $A$ into $k$ blocks $B_1,\dots,B_k$ such that $|s(B_i)-s(B_j| \le 1$ for every $i,j \in [k]$.
\end{theorem}

The {\bf proof} consists of checking that,  under conditions (i) and (ii), the algorithm for Theorem~\ref{th:main} works without any change.\hfill$\Box$

For instance, assume every $a_i$ is a $d$-dimensional vector with non-negative coordinates in the unit ball of the $\ell_p$ norm, $p\ge 1$. When $s(B)=||\sum_{a_i \in B}a_i||$, conditions (i) and (ii) are satisfied. (Simple examples show that this in not true for an arbitrary norm.) So there is a block partition $B_1,\dots,B_k$ of $A$ such that the norms of the sums of the elements in the blocks differ by at most one. However, unlike in the 1-dimensional case, this does not mean that corresponding vectors are
``almost'' equal.

Now we relax condition (ii):
\begin{enumerate}
\item[(iii)] $|s(B_1)- s(B_2)|\le 1$ if $B_1$ and $B_2$ differ by one element.
\end{enumerate}
In this case we can prove the same bound as in (\ref{eq:Mm}).

\begin{theorem}\label{th:gen} Assume $A=(a_1,\dots,a_n)$ and $k$ is positive integer. If $s$ is non-negative and satisfies conditions (i) and (iii), then there is a partition of $A$ into $k$ blocks $B_1,\dots,B_k$ such that $|s(B_i)-s(B_j)| \le 2$ for every $i,j \in [k]$.
\end{theorem}

Before the proof some preparation is in place. We are to consider intervals $[x,y)$ where $0\le x\le y\le n$. The interval $[i-1,i)$ is identified with the element $a_i$ of the sequence $A$.   The interval $[x,y)$ is a block if $x$ and $y$ are integers. Thus block $B=(a_i,\dots,a_j)$ can and will be identified with the interval $[i-1,j)$; here $i-1\le j$. Further, $[i,i)$ is the empty block positioned between $a_i$ and $a_{i+1}$.

Define next
\[
T^{k-1}_n=\{x=(x_1,\dots,x_{k-1})\in \R^{k-1}: 0\le x_1 \le \dots \le x_{k-1}\le n\},
\]
and set $x_0=0,x_k=n$. Every $x \in T_n^{k-1}$ determines a unique partition of $[0,n)$ into $k$ intervals
\[
[x_0,x_1),[x_1,x_2),\dots,[x_{k-1},x_k).
\]
Let $P(x)$ denote this partition. Also, conversely, every partition of $[0,n)$ into $k$ intervals determines a unique $x \in T_n^{k-1}$ such that $P(x)$ is equal to this partition. Note that $P(x)$ is a block partition if and only if all coordinates of $x$ are integers. In this case $P(x)=(B_1,\dots,B_k)$ and we define
\[
S(x)=(s(B_1),\dots,s(B_k)).
\]

The size $s(B)$ of block $B=(a_i,\dots,a_j)$ depends only on the interval $[i-1,j)$ so we may (and do) define $s([i-1,j))=s(B)$. For simpler notation we write $s[i,j)$ instead of $s([i,j))$. Note that $s[i,j)$ is always non-negative and $s[i,i)=0$. Condition (iii) says or rather implies that for all $0\le i\le j\le n$ and all $0\le i'\le j'\le n$
\[
|s[i,j)-s[i',j')| \le |i-i'|+|j-j'|.
\]
In other words, the map $s$ defined on pairs $(i,j)$ (with $0\le i\le j\le n$ ) is non-expanding in the $\ell_1$ norm.

We are going to extend $s$ from blocks $[i,j)$ to intervals $[x,y)$. Assume that $i,j \in [n]$, $i\le j$ and $x \in [i-1,i)$, $y \in [j-1,j)$.
The point $(x,y)\in \R^2$ is then either in the triangle with vertices $(i-1,j-1),(i,j-1),(i,j)$ or in the triangle with vertices $(i-1,j-1),(i-1,j),(i,j)$ or in both. Such triangles triangulate $T_n^2$ and so we can extend $s$ on each triangle linearly. This is the usual simplicial extension of $s$ onto $T_n^2$. We denote it invariably by $s$ so we have now an $s:T_n^2 \to \R$ map. It is very easy to check (and we omit the details) that the extended $s$ is also non-expanding, that is, for all $0\le x\le y\le n$ and all $0\le x'\le y'\le n$
\[
|s[x,y)-s[x',y')| \le |x-x'|+|y-y'|.
\]

The map $S$ was defined on the lattice points of $T_n^{k-1}$. We can extend it now to the whole $T_n^{k-1}$: for $x \in T_n^{k-1}$ let
\[
S(x)=(s[x_0,x_1),\dots,s[x_{k-1},x_k)) \in \R^k.
\]

{\bf Proof} of Theorem~\ref{th:gen}. Write $R(a,b)$ for the halfline starting at $a$ and going through $b$ where $a,b \in \R^k$ are distinct. Set $e=(1,1,\dots,1)\in \R^k$. We are going to show that there is an $x \in T_n^{k-1}$ such that $S(x)$ lies on the halfline $R(0,e)$. This is trivial if $S(x)$ coincides with the origin for some $x$. So we assume that $S(x)\ne 0$ for any $x \in T_n^{k-1}$.

Let $e_1,\dots,e_k$ be the standard basis of $R^k$ and write $\tri$ for the simplex with vertices $e_1,\dots,e_k$. We define a map $g: T_n^{k-1} \to \tri$ by setting $g(x)=R(0,S(x))\cap \tri$. As all coordinates of $S(x)$ are non-negative and $S(x)\ne 0$, $g$ is a continuous map.

The simplex $T_n^{k-1}$ has $k$ facets, $F_1,\dots,F_k$, where $F_i$ is given by the equation $x_{i-1}=x_i$. The facet $F_i$ is mapped by $g$ to points whose $i$th coordinate is zero. This implies that a $(k-1-h)$-dimensional face $F_{i_1}\cap \dots \cap F_{i_h}$ of $T_n^{k-1}$ is mapped by $g$ onto the $(k-1-h)$-face of $\tri$, defined by $z_{i_1}=\dots =z_{i_h}=0$ where $z_i$ is the $i$th coordinate of $z \in \tri$. In particular, $g$ is a one-to-one correspondence between the vertices of $T_n^{k-1}$ and the vertices of $\tri$. Let $f:\tri \to T_n^{k-1}$ be the linear (or simplicial) extension of $g^{-1}$ from the vertices of $\tri$ to the whole simplex $\tri$. It follows that $g\circ f:\tri \to \tri$ maps each face of $\tri$ onto itself.

\begin{lemma}\label{l:top} If $h: \tri \to \tri$ is continuous and maps each face of $\tri$ onto itself, then $h$ is surjective.
\end{lemma}

This result is known, see for instance Lemma 1 in \cite{kar} or Lemma 8.2 in \cite{kar1} and also \cite{Str}. For the convenience of the reader we give another short proof at the end of this paper.

The lemma implies that $g$ is also surjective. So there is an $x^* \in T_n^{k-1}$ with $g(x^*)=(1/k,\dots,1/k)$. Then $S(x^*)=(t,\dots,t)$ for some $t>0$ or for $t=0$ when $S(x^*)=0$ for some $x^* \in \tri$.

It is easy to finish the proof now. The point $x^*$ defines a partition $P(x^*)$ of $[0,n)$ into $k$ intervals $[x^*_{i-1},x^*_i)$ and $s[x^*_{i-1},x^*_i)=t$ for all $i \in [k]$. Round each $x_i^*$ to the nearest integer $y_i$, ties broken arbitrarily. So $|x_i^*-y_i|\le 1/2$. Now $y=(y_1,\dots,y_{k-1})$ defines a block partition $B_1,\dots,B_k$ of $[0,n)$ (or $A$, if you wish). As $s$ is non-expanding,
\begin{eqnarray*}
|s(B_i)&-&s(B_j)|=|s[y_{i-1},y_i)-s[y_{j-1},y_j)|\\
      &\le& |s[y_{i-1},y_i)-s[x^*_{i-1},x^*_i)|+|s[x^*_{i-1},x^*_i)-s[x^*_{j-1},x^*_j)|+\\
      &+&|s[x^*_{j-1},x^*_j)-s[y_{j-1},y_j)| \le 1+|t-t|+1 \le 2,
\end{eqnarray*}
for all $i,j \in [k]$, proving the theorem. \hfill$\Box$

\bigskip
{\bf Proof} of Lemma \ref{l:top}. Given such an $h$, the map $h_{\tau}:\tri \to \tri$ where $\tau \in [0,1]$ defined by $h_{\tau}(z)=(1-\tau)z+\tau h(z)$ is a homotopy between $h$ and the identity. Note that $h_{\tau}(z) \in \partial \tri$ if $z \in \partial \tri$ for all $\tau \in [0,1]$.

\medskip
Take now two disjoint copies, $\tri^+$ and $\tri^-$, of $\tri$ and identify their boundaries. This is an $S^d$, the $d$-dimensional unit sphere. Define a map $H:S^d \to S^d$ by setting $H(x)$ equal to $h(x)\in \tri^+$ if $x \in \tri^+$ and $h(x)\in \tri^-$ if $x \in \tri^-$. The map $H$ is well-defined and continuous since $x \in \partial \tri^{\pm}$ is mapped to $H(x) \in \partial \tri^{\pm}$. The homotopy $h_{\tau}$ extends to a homotopy $H_{\tau}$ between $H$ and the identity on $S^d$. Thus the degree of $H$ is one. This proves the lemma as $H|_{\tri^+}: \tri^+ \to \tri^+$ is the same as $h:\tri \to \tri$ and $H_{\tau}(z) \in \tri^+$ for all $\tau \in [0,1]$ if $z \in \tri^+$.
\hfill$\Box$

\bigskip
{\bf Acknowledgements.} We thank the anonymous referee for valuable comments. Research of the first author was partially supported by ERC Advanced Research Grant no 267165 (DISCONV), and by Hungarian National Research Grant K 83767.

\vspace{.5cm} {\sc Imre B\'ar\'any}\\
  {\footnotesize R\'enyi Institute of Mathematics}\\[-1.5mm]
  {\footnotesize Hungarian Academy of Sciences}\\[-1.5mm]
  {\footnotesize PO Box 127, 1364 Budapest, and }\\[-1.5mm]
  {\footnotesize Department of Mathematics}\\[-1.5mm]
  {\footnotesize University College London}\\[-1.5mm]
  {\footnotesize Gower Street, London WC1E 6BT}\\[-1.5mm]
  {\footnotesize England}\\[-1.5mm]
  {\footnotesize e-mail: {\tt barany@renyi.hu}}\\

{\sc Victor S. Grinberg}\\
{\footnotesize 5628 Hempstead Rd, Apt 102, Pittsburgh}\\[-1.5mm]
  {\footnotesize PA 15217, USA}\\[-1.5mm]
  {\footnotesize e-mail: {\tt victor}\_{\;}grinberg@yahoo.com}\\


\bigskip

\begin{thebibliography}{99}

\bibitem{BD} I. B\'ar\'any, B. Doerr. Balanced partitions of vector sequences, {\it Lin. Alg. Appl.}, {\bf 414} (2006), 464--469.

\bibitem{BG} I. B\'ar\'any, V. S. Grinberg. On some combinatorial questions in finite dimensional
spaces, {\it Lin. Alg.  Appl.} {\bf 41} (1981), 1--9.

\bibitem{kar} R. N. Karasev. KKM-type theorems for products of simplices and cutting sets and measures by straight lines (2009)
http://arxiv.org/pdf/0909.0604v1.pdf

\bibitem{kar1}R. N. Karasev. Geometry of measures: partitions and convex bodies, (2013)
http://www.rkarasev.ru/common/upload/mes\_partition.pdf

\bibitem{Str} N. Strickland. Map from simplex to itself that preserves subsimplices, (2011)
http://mathoverflow.net/questions/67318/

\end{thebibliography}
\end{document}